# SPHERICAL MAXIMAL OPERATORS ON RADIAL FUNCTIONS

Andreas Seeger, Stephen Wainger and James Wright

## 1. Introduction

For a function $f \in L^p(\mathbb{R}^d)$ we define the spherical means

$$A_t f(x) = \int_{S^{d-1}} f(x - ty') d\sigma(y')$$

where $d\sigma$ is the rotationally invariant measure on $S^{d-1}$, normalized such that $\sigma(S^{d-1}) = 1$.

Stein [5] showed that $\lim_{t \to 0} A_t f(x) = f(x)$ almost everywhere, provided $f \in L^p(\mathbb{R}^d)$, $p > d/(d-1)$ and $d \geq 3$. Later Bourgain [1] extended this result to the case $d = 2$. If $p \leq d/(d-1)$ then pointwise convergence fails. However if $\{t_j\}_{j=1}^\infty$ is a fixed sequence converging to 0 then pointwise convergence may hold for all $f \in L^p$ even if $p \leq d/(d-1)$, and $p$ depends on geometric properties of the sequence $\{t_j\}$.

According to a theorem by Stein [4] pointwise convergence holds for all $f \in L^p$ if the associated maximal operator

$$\mathcal{M}_E f(x) = \sup_{t \in E} |A_t f(x)|$$

is of weak type $(p, p)$, here $E = \{t_j : j \in \mathbb{N}\}$. Let $p(E)$ be the critical exponent for $L^p$-boundedness of $\mathcal{M}_E$, in the sense that $L^p$-boundedness holds for $p > p(E)$ and fails for $p < p(E)$. By the Marcinkiewicz interpolation theorem $p(E)$ is also the critical exponent for $M_E$ being of weak type $(p,p)$ if $1 < p \leq d/(d-1)$. A geometric characterization of $p(E)$ has been found in [3]; here arbitrary subsets $E$ of $(0, \infty)$ were admitted. In order to describe the result in [3] we let $I_k = [2^k, 2^{k+1}]$ and

$$E^k = I_k \cap E$$

and let $N(E^k, a)$ be the $a$-entropy number of $E^k$, that is the minimal number of intervals of length $a$ needed to cover $E^k$. Define

$$\kappa(E) = \overline{\lim_{\delta \to 0}} \Big[\sup_{k \in \mathbb{Z}} \frac{\log(1 + N(E^k, 2^k \delta))}{\log(1 + \delta^{-1})}\Big].$$

Then

$$p(E) = 1 + \frac{\kappa(E)}{d-1}.$$

Various results concerning the $L^p$-boundedness of $\mathcal{M}_E$ for the critical exponent $p = p(E)$ were proven in [3]; however these results fell short of being necessary and sufficient. If $p < d/(d-1)$ then a natural conjecture for the behavior on $L^p$ would be that $\mathcal{M}_E$ is of weak type $(p,p)$ if and only if the covering numbers $N(E^k, 2^k \delta)$ are bounded by $C\delta^{-(d-1)(p-1)}$, uniformly in $k$. Since the $L^p$-boundedness of $\mathcal{M}_E$ for $p < p(E)$ can be disproved by testing $\mathcal{M}_E$ on *radial* functions (in fact characteristic functions of balls) one might first examine the behavior of $\mathcal{M}_E$ on radial functions in $L^p$. In this paper we completely characterize the sets $E$ for which $\mathcal{M}_E$ is of strong type or of weak type $(p,p)$ on radial functions if $d \geq 3$ or if $d = 2$ and $p < 2$. Our first result concerns the case $p < d/(d-1)$.

The first and the second author were partially supported by an NSF grant; the third author was supported by an SERC grant.

**Theorem 1.1.** *Let $E \subset (0, \infty)$ and $d \geq 2$. Let $1 \leq p < \frac{d}{d-1}$. Then the inequality*

$$\left|\{x : \mathcal{M}_E f(x) > \alpha\}\right| \leq C \frac{\|f\|_{L^p_{\mathrm{rad}}}^p}{\alpha^p}$$

*hold for some $C$ and all radial $L^p$ functions $f$ if and only if for all $\delta \in (0, 1/2)$*

(1.1) $$\sup_{k \in \mathbb{Z}} N(E^k, 2^k \delta) \leq C' \delta^{-(d-1)(p-1)}$$

*for some $C'$ independent of $\delta$.*

The condition for a strong type inequality is somewhat more complicated. More generally we consider the $L^p_{\mathrm{rad}} \to L^{pq}$ mapping properties where $L^{pq}$ is the standard Lorentz-space.

**Theorem 1.2.** *Let $E \subset (0, \infty)$, $d \geq 2$ and let $1 \leq p < d/(d-1)$, $p \leq q < \infty$. Then the inequality*

$$\|\mathcal{M}_E f\|_{L^{pq}(\mathbb{R}^d)} \leq C \|f\|_{L^p_{\mathrm{rad}}(\mathbb{R}^d)}$$

*holds for some $C$ and all radial $f \in L^p(\mathbb{R}^d)$ if and only if the condition*

(1.2) $$\sup_{j \in \mathbb{Z}} \Big(\sum_{n=0}^{\infty} [N(E^{j+n}, 2^j)]^{q/p} 2^{-n(d-1)q/p'}\Big)^{1/q} < \infty$$

*is satisfied.*

Note that Theorem 1.1 is just the limiting case of Theorem 1.2, for $q \to \infty$.

In dimensions $d \geq 3$ we can also prove a characterization for $L^p \to L^{pq}$-boundedness at the critical exponent $p = (d-1)/d$, on radial functions. For the case $q > p$ we have

**Theorem 1.3.** *Let $E \subset (0, \infty)$ and $d \geq 3$. Let $p_d = d/(d-1)$. Then $\mathcal{M}_E$ is of weak type $(p_d, p_d)$ on radial functions if and only if for all $\delta \in (0, 1/2)$*

(1.3) $$\sup_{k \in \mathbb{Z}} N(E^k, 2^k \delta) \leq C \delta^{-1} \big[\log(1/\delta)\big]^{-1/(d-1)}.$$

If $p_d < q < \infty$ then $\mathcal{M}_E$ maps $L^{p_d}_{\mathrm{rad}}$ boundedly into $L^{p_d q}$ if and only if for all $\delta \in (0, 1/2)$

(1.4) $$\sup_{|I| = \log \delta^{-1}} \Big(\sum_{k \in I} [N(E^k, 2^k \delta)]^{q/p_d}\Big)^{p_d/q} \leq C \delta^{-1} \big[\log(1/\delta)\big]^{-1/(d-1)}$$

where the supremum is taken over all intervals $I$ of length $\log \delta^{-1}$.

The condition for $q = p_d$ takes a different form. Let $\delta_{(k,n)}$ be the Dirac measure in $\mathbb{R}^2$ supported at $(k, n)$. For any subinterval $I$ of the real line let $T(I)$ be the tent

(1.5) $$T(I) = \{(x, t) : x \in I, \, 0 \leq t \leq |I|\}.$$

**Theorem 1.4.** *If $d \geq 3$ then $\mathcal{M}_E$ maps $L^{pq}_{\text{rad}}(\mathbb{R}^d)$ boundedly into $L^{pq}(\mathbb{R}^d)$ if and only if the discrete measure $\sum_{k \in \mathbb{Z}} \sum_{n>0} N(E^k, 2^{k-n}) 2^{-n} n^{1/(d-1)} \delta_{(k,n)}$ on the upper half plane is a Carleson measure; i.e.*

$$(1.6) \qquad \sup_{|I| \geq 1} \frac{1}{|I|} \sum_{(k,n) \in T(I)} N(E^k, 2^{k-n}) 2^{-n} n^{1/(d-1)} < \infty.$$

For a different formulation of Theorems 1.3 and 1.4 unifying the statements for $p \leq q \leq \infty$ see Corollary 2.6 below.

*Remarks.* In various instances the conditions for boundedness take a simpler form, see also the discussion in §2.

(i) Let $E = E^0$ be supported in $[1,2]$. Let

$$(1.7) \qquad w_n = |\{r : \text{dist}(r, E^0) \leq 2^{-n}\}|.$$

Then for $1 \leq p < d/(d-1)$ the local maximal operator $\mathcal{M}_{E^0}$ is bounded on $L^p_{\text{rad}}$ if and only if the inequality

$$(1.8) \qquad w_n \leq C 2^{-n(1-(d-1)(p-1))}$$

holds uniformly in $n \geq 1$. Moreover $\mathcal{M}_{E^0}$ is bounded on $L^{d/(d-1)}_{\text{rad}}$ if and only if the inequality

$$(1.9) \qquad w_n \leq C n^{-\frac{1}{d-1}}$$

holds uniformly in $n \geq 1$. These facts are immediate consequences of Theorems 1.2 and 1.4. Note also that in this case the condition for $L^p_{\text{rad}} \to L^{pq}$-boundedness does not depend on $q \in [p, \infty]$.

(ii) Suppose that $E$ is the union of $2^k$-dilates of a fixed set $E^0$ supported in $[1,2]$ and let $w_n$ be defined as in (1.7). Then the condition

$$(1.10) \qquad \left( \sum_{n=0}^{\infty} \left[ w_n 2^{n(1-(d-1)(p-1))} \right]^{q/p} \right)^{1/q} < \infty$$

is equivalent to (1.2) and the condition

$$(1.11) \qquad w_n \leq C n^{-\frac{1}{d-1}(1+\frac{d}{q})}$$

is equivalent to (1.4) (or (1.3) respectively). Finally condition (1.6) becomes

$$(1.12) \qquad \sum_{n=1}^{\infty} w_n n^{1/(d-1)} < \infty.$$

Various equivalent forms of our conditions are discussed in §2. There the necessity of these conditions is also proved. The proofs of the sufficiency are contained in §3-5. The proofs of Theorem 1.1 and 1.2 are contained in §3 for the case $d = 3$ and in §5 for the case $d = 2$. The proofs of Theorems 1.3 and 1.4 are given in §4.

The following notation is used: For a set $E \subset \mathbb{R}^d$ we denote the Lebesgue measure of $E$ by $|E|$. The measure $\mu_d$ on $\mathbb{R}_+$ is defined by $d\mu_d = r^{d-1} dr$. We always let $I_k = [2^k, 2^{k+1}]$ and $\widetilde{I}_k = [2^{k-1}, 2^{k+2}]$. Given two quantities $a$ and $b$ we write $a \lesssim b$ or $b \gtrsim a$ if there is a positive constant $C$ such that $a \leq Cb$. We write $a \approx b$ if $a \lesssim b$ and $a \gtrsim b$.

## 2. Preliminaries

We begin by recalling a characterization of Lorentz-spaces $L^{pq}$. Let $f$ be a measurable function in a measure space $\Omega$ with measure $\mu$. Let $\lambda_f(\alpha) = \mu(\{x : |f(x)| > \alpha\})$ be the distribution function and $f^*(t) = \inf\{\alpha : \lambda_f(\alpha) \leq t\}$ the nonincreasing rearrangement. Then one defines the $L^{pq}$-quasinorm with respect to the measure $\mu$ by

$$\|f\|_{L^{pq}(\mu)} = \left(\frac{q}{p}\int_0^\infty \left[t^{1/p}f^*(t)\right]^q \frac{dt}{t}\right)^{1/q}.$$

This is not actually a norm but for $1 < p < \infty$ the space $L^{pq} = \{f : \|f\|_{L^{pq}} < \infty\}$ carries a Banach space topology ([7]). The following Lemma is probably well known, but we include a proof because of lack of an appropriate reference.

**Lemma 2.1.** *For any measurable function $f$*

$$(2.1) \qquad \|f\|_{L^{pq}} = \left(q\int_0^\infty \alpha^q [\lambda_f(\alpha)]^{q/p} \frac{d\alpha}{\alpha}\right)^{1/q} \approx \left(\sum_{\sigma \in \mathbb{Z}} 2^{\sigma q}[\lambda_f(2^\sigma)]^{q/p}\right)^{1/q}$$

**Proof.** The equivalence of the second and third expression follows from a standard argument (as in [7, p.192]). It suffices to show the equality of the first two terms for nonnegative simple functions; the general case follows by a limiting argument (see [7, p.191]). Therefore assume $f(x) = \sum_{j=1}^n a_j \chi_{E_j}(x)$, with $a_1 \geq a_2 \cdots \geq a_n \geq 0$. Let $v_j = \sum_{k=1}^j \mu(E_k)$. Set $v_0 = 0$ and $a_{n+1} = 0$. Then $\lambda_f(\alpha) = \sum_{j=1}^n v_j \chi_{[a_{j+1},a_j)}(\alpha)$ and $f^*(t) = \sum_{j=1}^n a_j \chi_{[v_{j-1},v_j)}(t)$. An evaluation of the integrals yields

$$\frac{q}{p}\int_0^\infty [t^{1/p}f^*(t)]^q \frac{dt}{t} = \sum_{j=1}^n a_j^q [v_j^{q/p} - v_{j-1}^{q/p}]$$

$$q\int_0^\infty \alpha^q [\lambda_f(\alpha)]^{q/p}\frac{d\alpha}{\alpha} = \sum_{j=1}^n v_j^{q/p}[a_j^q - a_{j+1}^q]$$

and since $v_0 = 0$, $a_{n+1} = 0$ the two expressions coincide. $\square$

It is sometimes useful to express the conditions in Theorems 1.1 and 1.2 in different ways. Let $N(E^k, 2^{k-n})$ be as in the introduction and define

$$(2.2) \qquad W_n^k = \{r \in I_k : \text{dist}(r, E^k) \leq 2^{k-n+1}\}$$

and

$$(2.3) \qquad D_n^k = \{r \in I_k : 2^{k-n} < \text{dist}(r, E^k) \leq 2^{k-n+1}\}.$$

By a binary interval of length $2^j$ we mean an interval of the form $[m2^j, (m+1)2^j]$ for $m, j \in \mathbb{Z}$. We first note

**Lemma 2.2.** *Let $E \subset (0,\infty)$. Let $\widetilde{N}(E, 2^j)$ be the minimal number of binary intervals of length $2^j$ needed to cover $E$. Then*

$$N(E, 2^j) \leq \widetilde{N}(E, 2^j) \leq 3N(E, 2^j).$$

*Moreover for any interval $I$ with $2^j \leq |I|$*

$$2^{j-2}N(E \cap I, 2^j) \leq \left|\{r \in I : \text{dist}(r, E \cap I) \leq 2^{j+1}\}\right| \leq 2^{j+3}N(E \cap I, 2^j).$$

We omit the elementary proof.

**Lemma 2.3.** *Let $E \subset (0, \infty)$ such that $|\overline{E}| = 0$. Suppose that $1 \leq p < d/(d-1)$ and $q \geq 1$. Then the following conditions are equivalent:*

$$\sup_{j \in \mathbb{Z}} \Big( \sum_{n=0}^{\infty} \Big[ N(E^{n+j}, 2^j) \Big]^{1/p} 2^{-n(d-1)/p'} \Big]^q \Big)^{1/q} < \infty \tag{2.4}$$

$$\sup_{j \in \mathbb{Z}} 2^{-j/p} \Big( \sum_{n=0}^{\infty} \Big[ |W_n^{j+n}|^{1/p} 2^{-n(d-1)/p'} \Big]^q \Big)^{1/q} < \infty \tag{2.5}$$

$$\sup_{j \in \mathbb{Z}} 2^{-j/p} \Big( \sum_{n=0}^{\infty} \Big[ |D_n^{j+n}|^{1/p} 2^{-n(d-1)/p'} \Big]^q \Big)^{1/q} < \infty. \tag{2.6}$$

*For $q = \infty$ these statements remain true if one replaces the $\ell^q$-norm by the supremum.*

**Proof.** The equivalence of (2.4) and (2.5) for all $p$ immediately follows from Lemma 2.2. Clearly (2.5) implies (2.6). Since we assume that the closure of $E$ is a null set we also have

$$|W_n^k| = \sum_{m \geq n} |D_m^k|. \tag{2.7}$$

Using (2.7) and Minkowski's inequality we obtain

$$2^{-j/p} \Big( \sum_{n=0}^{\infty} \Big[ |W_n^{j+n}|^{1/p} 2^{-n(d-1)/p'} \Big]^q \Big)^{1/q}$$

$$\leq 2^{-j/p} \sum_{s \geq 0} \Big( \sum_{n=0}^{\infty} \Big[ |D_{n+s}^{j+n}|^{1/p} 2^{-n(d-1)/p'} \Big]^q \Big)^{1/q}$$

$$\leq \sum_{s \geq 0} 2^{s[(d-1)/p' - 1/p]} \sup_{l \in \mathbb{Z}} \Big[ 2^{-(l-s)/p} \Big( \sum_{m=s}^{\infty} |D_m^{l-s+m}|^{q/p} 2^{-m(d-1)q/p'} \Big)^{1/q} \Big].$$

Since $j$ was arbitrary and $p < d/(d-1)$ we see that (2.6) implies (2.5). $\square$

We now discuss alternative formulations of the conditions in Theorems 1.3 and 1.4. Define for any subset $A$ of $\mathbb{Z} \times \mathbb{N}$

$$\mathfrak{v}_\lambda(A) = \sum_{(k,n) \in A} |D_n^k| 2^{k(\lambda-1)} \tag{2.8}$$

where in our application $\lambda = d$. Let

$$\Gamma_\lambda(\beta) = \{(k,n) : n 2^{-k(\lambda-1)} > \beta\} \tag{2.9}$$

and for any interval $I$ let $T(I)$ be the tent of $I$ as defined in (1.5).

**Lemma 2.4.** *Suppose $E \subset (0, \infty)$ and $|\overline{E}| = 0$. Suppose $1 < p < \infty$, $1 < \lambda < \infty$. Then*

$$\sup_{n \geq 1} \sup_{k \in \mathbb{Z}} |W_n^k| 2^{-k(\lambda-1)(p-1)} n^{p-1} \approx \sup_{|I| \geq 1} \sup_\beta \beta^p \frac{\mathfrak{v}_\lambda(\Gamma_\lambda(\beta) \cap T(I))}{|I|}$$

*and for $1 < p < q < \infty$*

$$\sup_{n \geq 1} \sup_{|I|=n} \Big( \sum_{k \in I} \big[ |W_n^k| 2^{-k(\lambda-1)(p-1)} n^{p-1} \big]^{q/p} \Big)^{p/q} \approx \sup_{|I| \geq 1} \Big( \int_0^\infty \beta^q \Big[ \frac{\mathfrak{v}_\lambda(\Gamma_\lambda(\beta) \cap T(I))}{|I|} \Big]^{q/p} \frac{d\beta}{\beta} \Big)^{p/q}.$$

**Proof.** We only consider the case $p < q < \infty$, the case $q = \infty$ is proved in the same way.

First fix $n \geq 1$ and an interval $I = [a, b]$, $b - a = n$. For $l \geq 1$ let $I^l = [a - 2^l n, b + 2^l n]$ and let

$$\beta_l(k, n) = 2^{l-1} n 2^{-k(\lambda - 1)}$$

so that

(2.10) $$\frac{\beta_l(k, n)}{\beta_l(k+1, n)} = 2^{\lambda - 1} \neq 1.$$

Then using (2.7) we obtain

$$\Big(\sum_{k \in I} [2^{-k(\lambda-1)(p-1)}|W_n^k|n^{p-1}]^{q/p}\Big)^{p/q}$$

$$\leq n^{p-1} \sum_{l=1}^{\infty} \Big(\sum_{k \in I} [2^{-k(\lambda-1)p} \sum_{m=2^{l-1}n}^{2^l n} |D_m^k| 2^{k(\lambda-1)}]^{q/p}\Big)^{p/q}$$

$$\leq n^{p-1} \sum_{l=1}^{\infty} \Big(\sum_{k \in I} [(\beta_l(k,n) 2^{1-l}/n)^p \mathfrak{v}_\lambda(T(I_l) \cap \Gamma_\lambda(\beta_l(k,n)))]^{q/p}\Big)^{p/q}$$

$$\lesssim \sum_{l=1}^{\infty} 2^{-l(p-1)} \Big(\sum_k \Big[(\beta_l(k,n))^p \frac{\mathfrak{v}_\lambda(T(I_l) \cap \Gamma_\lambda(\beta_l(k,n)))}{|I_l|}\Big]^{q/p}\Big)^{p/q}$$

$$\lesssim \sup_{|I| \geq 1} \Big(\int_0^\infty \beta^q \Big[\frac{\mathfrak{v}_\lambda(T(I) \cap \Gamma_\lambda(\beta))}{|I|}\Big]^{q/p} \frac{d\beta}{\beta}\Big)^{p/q}$$

where for the last inequality we have used (2.10).

For the opposite inequality we fix an interval $I$ of length $|I| \geq 1$ and set for $\beta > 0$, $m \geq 1$, $1 \leq 2^s \leq 2|I|$

$$\mathfrak{A}^\beta_{sm} = \{k \in I : \beta 2^{m-2-s} \leq 2^{-k(\lambda-1)} < \beta 2^{m-s+1}\}$$

and observe that the cardinality of $\mathfrak{A}^\beta_{sm}$ is bounded independently of $s$, $m$ and $\beta$. Moreover for fixed $k$, $s$ and $m$

(2.11) $$\int_{\{\beta : k \in \mathfrak{A}^\beta_{sm}\}} \frac{d\beta}{\beta} \leq C_\lambda.$$

We estimate

$$\Big(\int_0^\infty \beta^q \Big[|I|^{-1} \sum_{(k,n) \in T(I) \cap \Gamma_\lambda(\beta)} 2^{k(\lambda-1)}|D_n^k|\Big]^{q/p} \frac{d\beta}{\beta}\Big)^{p/q}$$

$$\lesssim \sum_{m=0}^\infty \sum_{s=1}^{[1+\log_2 |I|]} \Big(\int_0^\infty \beta^q \Big[|I|^{-1} \sum_{k \in \mathfrak{A}^\beta_{sm}} \sum_{n=2^{s-1}}^{2^s} 2^{-k(\lambda-1)(p-1)}|D_n^k|(2^{s-m}/\beta)^p\Big]^{q/p} \frac{d\beta}{\beta}\Big)^{p/q}$$

$$\lesssim \sum_{m=0}^\infty 2^{-mp} \sum_{s=1}^{[1+\log_2 |I|]} 2^s |I|^{-1} \Big(\int_0^\infty \Big[\sum_{k \in \mathfrak{A}^\beta_{sm}} 2^{-k(\lambda-1)(p-1)}|W^k_{2^{s-1}}|2^{s(p-1)}\Big]^{q/p} \frac{d\beta}{\beta}\Big)^{p/q}$$

$$\lesssim \sum_{s=1}^{[1+\log_2 |I|]} 2^s |I|^{-1} \Big(\sum_{k \in I} [2^{-k(\lambda-1)(p-1)}|W^k_{2^{s-1}}|2^{s(p-1)}]^{q/p}\Big)^{p/q}$$

where for the last inequality we have used (2.11) and the statement preceding it. Now for each $s$ with $1 \leq 2^s \leq 2|I|$ the interval $I$ is a union of $\approx |I|2^{-s}$ intervals $J_\nu^s$ of length $2^s$ and therefore

$$\sum_{s=1}^{[1+\log_2 |I|]} 2^s |I|^{-1} \Big(\sum_{k \in I} [2^{-k(\lambda-1)(p-1)}|W_{2^{s-1}}^k|2^{s(p-1)}]^{q/p}\Big)^{p/q}$$

$$\lesssim \sum_{s=1}^{[1+\log_2 |I|]} 2^s |I|^{-1} \Big(\sum_\nu \sum_{k \in J_\nu^s} [2^{-k(\lambda-1)(p-1)}|W_{2^{s-1}}^k|2^{s(p-1)}]^{q/p}\Big)^{p/q}$$

$$\lesssim \sum_{s=1}^{[1+\log_2 |I|]} 2^{s(1-p/q)} |I|^{-1+p/q} \sup_{n \geq 1} \sup_{|J|=n} \Big(\sum_{k \in J} [2^{-k(\lambda-1)(p-1)}|W_n^k|n^{p-1}]^{q/p}\Big)^{p/q}$$

and the desired inequality follows since $q > p$. $\square$

**Lemma 2.5.** *Suppose $1 < p < \infty$, $0 < \lambda < \infty$. Then the condition*

$$(2.12) \qquad \sup_{|I| \geq 1} \frac{1}{|I|} \sum_{(k,n) \in T(I)} |W_n^k| 2^{-k(\lambda-1)(p-1)} n^{p-1} < \infty$$

*holds if and only if $|\overline{E}| = 0$ and the condition*

$$(2.13) \qquad \sup_{|I| \geq 1} \int_0^\infty \beta^p \frac{\mathfrak{v}_\lambda(\Gamma_\lambda(\beta) \cap T(I))}{|I|} \frac{d\beta}{\beta} < \infty$$

*holds.*

**Proof.** We first observe that

$$p \int_0^\infty \beta^p \frac{\mathfrak{v}_\lambda(\Gamma_\lambda(\beta) \cap T(I))}{|I|} \frac{d\beta}{\beta} = p \frac{1}{|I|} \int_0^\infty \beta^p \sum_{(k,n) \in \Gamma_\lambda(\beta) \cap T(I)} |D_n^k| 2^{k(\lambda-1)} \frac{d\beta}{\beta}$$

$$= \frac{1}{|I|} \sum_{(k,n) \in T(I)} |D_n^k| 2^{k(\lambda-1)} \int_0^{n 2^{-k(\lambda-1)}} p \beta^{p-1} d\beta$$

$$(2.14) \qquad = \frac{1}{|I|} \sum_{(k,n) \in T(I)} |D_n^k| 2^{-k(\lambda-1)(p-1)} n^p.$$

Now suppose that (2.12) holds. Then it is easy to see that $|\overline{E}| = 0$. Fix an interval $I$ and let $I^*$ the double interval. Since the sequence $n \to |W_n^k|$ is monotone we obtain that

$$\frac{1}{|I|} \sum_{(k,n) \in T(I)} |D_n^k| 2^{-k(\lambda-1)(p-1)} n^p \lesssim \frac{1}{|I|} \sum_{k \in I} \sum_{s=0}^{[1+\log_2(|I|)]} |W_{2^s}^k| 2^{-k(\lambda-1)(p-1)} 2^{sp}$$

$$(2.15) \qquad \lesssim \frac{1}{|I^*|} \sum_{(k,n) \in T(I^*)} |W_n^k| 2^{-k(\lambda-1)(p-1)} n^{(p-1)}.$$

Conversely assume that $|\overline{E}| = 0$ so that (2.7) applies. Fix $I = [a, b]$ and let $I_l = [a - 2^{l+1}|I|, b + 2^{l+1}|I|]$.

Then

$$\frac{1}{|I|} \sum_{(k,n) \in T(I)} |W_n^k| 2^{-k(\lambda-1)(p-1)} n^{p-1}$$

$$\lesssim \frac{1}{|I|} \sum_{k \in I} \sum_{s=0}^{[\log_2(|I|)+1]} 2^{-k(\lambda-1)(p-1)} |W_{2^s}^k| 2^{sp}$$

$$\lesssim \frac{1}{|I|} \sum_{l=0}^{\infty} \sum_{k \in I} \sum_{s=0}^{[\log_2(|I|)+1]} 2^{-k(\lambda-1)(p-1)} 2^{sp} \sum_{2^{s+l-1} \leq m \leq 2^{s+l}} |D_m^k|$$

$$\lesssim \sum_{l=0}^{\infty} \frac{2^l}{|I_l|} 2^{-lp} \sum_{k \in I} \sum_{s=0}^{[\log_2(|I|)+1]} 2^{-k(\lambda-1)(p-1)} \sum_{2^{s+l-1} \leq m \leq 2^{s+l}} |D_m^k| m^p$$

$$\lesssim \sum_{l=0}^{\infty} 2^{-l(p-1)} \frac{1}{|I_l|} \sum_{(k,m) \in T(I_l)} 2^{-k(\lambda-1)(p-1)} |D_m^k| m^p$$

(2.16) $$\lesssim \sup_J \frac{1}{|J|} \sum_{(k,m) \in T(J)} 2^{-k(\lambda-1)(p-1)} |D_m^k| m^p.$$

The asserted equivalence follows from (2.14-16). □

Lemmas 2.2, 2.4 and 2.5 can be used to unify the statements of Theorems 1.3 and 1.4.

**Corollary 2.6.** *For any interval $I$ with $|I| \geq 1$ let*

(2.17) $$\mathfrak{a}_{p,I}(\beta) = \beta \Big[\frac{\mathfrak{v}_d(\Gamma_d(\beta) \cap T(I))}{|I|}\Big]^{1/p}.$$

*Then the condition*

(2.18) $$\sup_{|I| \geq 1} \|\mathfrak{a}_{p_d,I}\|_{L^q(\mathbb{R}_+, d\beta/\beta)} < \infty$$

*is equivalent with (1.3) if $q = \infty$, with (1.4) if $p < q < \infty$ and with (1.6) if $q = p$.*

**Necessary conditions.** We now consider the spherical mean $A_t f$ for a *radial* function $f$ with $f(x) = f_0(|x|)$. Then $A_t f$ is also a radial function, given by

(2.19) $$A_t f(x) = c_d \int_{||x|-t|}^{|x|+t} K_t(|x|, s) f_0(s) ds$$

where

(2.20) $$K_t(r,s) = \Big[\frac{\sqrt{(r+t)^2 - s^2}\sqrt{s^2 - (r-t)^2}}{(r+t)^2 - (r-t)^2}\Big]^{d-3} \frac{s}{(r+t)^2 - (r-t)^2}.$$

This follows from a straightforward computation, see [2]. In order to derive necessary conditions for $L^p \to L^{pq}$ boundedness we shall use the following lower bounds which immediately follow from (2.19-20).

**Lemma 2.7.** *Suppose $f_0(r) \geq 0$ for all $r > 0$. Then there is $c > 0$, independent of $f_0$, such that*

$$\mathcal{M}_E f(x) \geq c 2^{-k(d-1)} \int_{2^{k-n+2}}^{2^k} s^{d-2} f_0(s) ds \quad \text{if } |x| \in D_n^k.$$

We first note set if $f_0(s) = s^{-(d-1)}[\log 1/s]^{-1}\chi_{[0,1/2]}(s)$ then $f \in L^p(\mathbb{R}^d)$, $p \leq d/(d-1)$, but $\mathcal{M}_E f(x) = \infty$ for $x \in \overline{E}$. Thus $L^p \to L^{pq}$ boundedness for $p \leq d/(d-1)$ implies that $\overline{E}$ is a null set.

In order to see the sharpness of Theorem 1.1 and Theorem 1.2 we test $\mathcal{M}_E$ on the radial $f$ with

$$f_0(s) = s^{-d/p}\chi_{[2^j, 2^{j+1}]}(s),$$

then $\|f\|_p \approx 1$. Also

$$|\{x : \mathcal{M}_E f(x) > C 2^\sigma\}| \gtrsim \sum_{n : 2^{-n(d-1)}2^{-jd/p} > 2^\sigma} |D_n^{j+n}| 2^{(j+n)(d-1)}$$

for suitable $C > 0$. We estimate from below the sum of the right hand side by the sum over those terms with $2^{-n(d-1)}2^{-jd/p} \approx 2^{\sigma+10}$ and use only those expressions in the definition of the Lorentz-space via the distribution function. This yields

$$\left(\sum_\sigma [2^\sigma |\{x : \mathcal{M}_E f(x) > C 2^\sigma\}|^{1/p}]^q\right)^{1/q} \gtrsim \left(\sum_{n=0}^\infty [|D_n^{j+n}|^{1/p} 2^{-j/p} 2^{-n(d-1)/p'}]^q\right)^{1/q}.$$

Since $j$ is arbitrary the necessity of the conditions in Theorems 1.1 and 1.2 follows from the last inequality and Lemma 2.3.

In order to see the sharpness of Theorems 1.3 and 1.4 we fix $I = [a - L, a + L]$ and define $f$ via

$$f_0(s) = s^{1-d}\chi_{[2^{a-10L}, 2^{a+10L}]}(s).$$

Then

(2.21) $$\|f\|_{L^{p_d}(\mathbb{R}^d)} = O(|I|^{(d-1)/d}).$$

Now if $|x| \in D_n^k$, $k \in I$, $0 < n \leq |I|$ then $\mathcal{M}_E f(x) \geq c 2^{-k(d-1)} n$ by Lemma 2.7 and therefore

(2.22) $$|\{x \in \mathbb{R}^d : \mathcal{M}_E f(x) > \alpha\}| \geq \mathfrak{v}_d(\Gamma_d(c\alpha) \cap T(I)).$$

In view of Corollary 2.6 the necessity of the conditions in Theorem 1.3 and Theorem 1.4 is an immediate consequence of (2.21) and (2.22).

## 3. Estimates in higher dimensions

We shall use the following pointwise estimate for the spherical means acting on radial functions $f$ defined in $\mathbb{R}^d$, $d \geq 3$ such that

$$f(x) = f_0(r) \text{ where } r = |x|.$$

**Lemma 3.1.** *Fix $1 \leq p < 2$ and set*

(3.1) $$g(s) = f_0(s)s^{(d-1)/p}.$$

*Then*

(3.2) $$\mathcal{M}_E f(x) \leq C_1[\mathfrak{M}g(r) + R_1 f_0(r) + R_2 f_0(r)]$$

*where*

(3.3) $$\mathfrak{M}g(r) = \sup_{\substack{t \in E \\ r/2 < t < 3r/2}} r^{1-d} \int_{|r-t|}^{r+t} s^{\frac{d-1}{p'}-1} g(s) ds$$

(3.4) $$R_1 f_0(r) = \sup_{t \leq r/2} \frac{1}{t} \int_{r-t}^{r+t} f_0(s) ds$$

(3.5) $$R_2 f_0(r) = \sup_{t \geq 3r/2} \frac{1}{r} \int_{t-r}^{t+r} f_0(s) ds.$$

The estimate (3.2) is an easy consequence of (2.19-20); we shall omit the proof. Theorems 1.1 and 1.2 for $d \geq 3$ are immediate consequences of the estimates for $R_1$, $R_2$ and $\mathfrak{M}$ in Lemmas 3.2, 3.3 and Proposition 3.4 below.

**Lemma 3.2.** *For $1 < p \leq \infty$, $d \geq 1$ the operator $R_1$ is bounded on $L^p(\mu_d)$. Moreover $R_1$ is of weak type $(1,1)$, with respect to the measure $\mu_d$; i.e.*

$$\mu_d(\{r : R_1 f_0(r) > \alpha\}) \lesssim \alpha^{-1} \int |f_0(s)| s^{d-1} ds,$$

*for all $\alpha > 0$.*

**Proof.** Since $R_1$ is bounded on $L^\infty$ it suffices to prove the weak-type inequality and the conclusion follows by the Marcinkiewicz interpolation theorem. We observe that if $\operatorname{supp} h \subset [2^k, 2^{k+1}]$ then $\operatorname{supp} R_1 h \subset [2^{k-1}, 2^{k+2}]$. Let $f^k(s) = \sum_{i=-2}^{2} |f_0(s)| \chi_{I_{k+i}}(s)$. Then

$$\mu_d(\{r : R_1 f_0(r) > \alpha\}) \lesssim \sum_k 2^{k(d-1)} |\{r \in I_k : R_1 f^k(r) > \alpha\}|$$

and by the weak type inequality for the Hardy-Littlewood maximal function the right hand side is dominated by a constant times

$$\alpha^{-1} \sum_k 2^{k(d-1)} \|f^k\|_1 \lesssim \alpha^{-1} \int |f_0(s)| s^{d-1} ds. \quad \square$$

**Lemma 3.3.** *Suppose $d > 1$. Then the operator $R_2$ is bounded on $L^p(\mu_d)$, $1 \leq p \leq \infty$.*

**Proof.** Clearly $R_2$ is bounded on $L^\infty$, so it suffices to check the claim for $p = 1$. Then

$$\sum_{k \in \mathbb{Z}} \int_{2^k}^{2^{k+1}} |R_2 f_0(r)| r^{d-1} dr \lesssim \sum_{k \in \mathbb{Z}} \sum_{L \geq 0} \int_{2^k}^{2^{k+1}} r^{d-2} 2^{-(k+L)(d-1)} \int_{2^{k+L-1}}^{2^{k+L+1}} |f_0(s)| s^{d-1} ds\, dr$$

$$\lesssim \sum_{L \geq 0} 2^{-L(d-1)} \sum_{k \in \mathbb{Z}} \int_{2^{k+L-1}}^{2^{k+L+1}} |f_0(s)| s^{d-1} ds \lesssim \|f_0\|_{L^1(\mu_d)}. \quad \square$$

The assumptions on the set $E$ will be needed now when we estimate $\mathfrak{M} f_0$.

**Proposition 3.4.** *Suppose that $1 \le p < d/(d-1)$, $p \le q \le \infty$ and*

$$\sup_j \Big(\sum_{n \ge 0}\big[N(E^{j+n}, 2^j)^{1/p} 2^{-n(d-1)/p'}\big]^q\Big)^{1/q} \le A \qquad \text{if } p \le q < \infty$$

$$\sup_{j,n} N(E^{j+n}, 2^j)^{1/p} 2^{-n(d-1)/p'} \le A \qquad \text{if } q = \infty.$$

*Then $\mathfrak{M}$ maps $L^p$ boundedly into $L^{pq}(\mu_d)$.*

**Proof.** Define for $\ell \ge 0$

$$\mathfrak{M}_\ell g(r) = \sum_k \sum_{n \ge \ell-3} \chi_{D_n^k}(r) r^{1-d} \int_{2^{k-n+\ell}}^{2^{k-n+\ell+1}} s^{\frac{d-1}{p'}-1} g(s)\,ds; \tag{3.6}$$

then $\mathfrak{M}g(r) \le \sum_{\ell=0}^\infty \mathfrak{M}_\ell g(r)$. We now derive an $L^{pq}(\mu_d)$-estimate for $\mathfrak{M}_\ell g$ in terms of the $L^p$-norm of $g$ (which is equal to the $L^p(\mu_d)$-norm of $f_0$).

First note that the assumption on $E$ implies $|\overline{E}| = 0$ and therefore also $\mu_d(\overline{E}) = 0$. Consequently it suffices to estimate the functions $\mathfrak{M}_\ell g$ on the set $\cup_{k,n} D_n^k$. By Hölder's inequality

$$\mathfrak{M}_\ell g(r) \le C_1 2^{-k(d-1)} 2^{(k-n+\ell)(\frac{d-1}{p'}-\frac{1}{p})} \|g\|_{L^p(I^{k-n+\ell})} \qquad \text{if } r \in D_n^k.$$

Therefore

$$\mu_d(\{r : \mathfrak{M}_\ell g > 2^\sigma\}) \lesssim \sum_* 2^{k(d-1)} |D_n^k|$$

where in the starred sum we sum over all pairs $(k,n)$ with the property that

$$C_1 2^{-k(d-1)} 2^{(k-n+\ell)(\frac{d-1}{p'}-\frac{1}{p})} \|g\|_{L^p(I^{k-n+\ell})} > 2^\sigma.$$

We change variables $j = k - n$ and set for fixed $j$ and $\kappa = 0, 1, 2, \dots$

$$\mathcal{B}_{\ell j \sigma}^\kappa = \{n \ge 0 : C_1 2^{-(j+n)(d-1)} 2^{(j+\ell)(\frac{d-1}{p'}-\frac{1}{p})} \|g\|_{L^p(I^{j+\ell})} \in [2^{\sigma+\kappa}, 2^{\sigma+\kappa+1})\}$$

$$= \{n \ge 0 : C_1^p 2^{-jd} 2^{-n(d-1)p} 2^{\ell(\frac{d-1}{p'}-\frac{1}{p})p} \|g\|_{L^p(I^{j+\ell})}^p \in [2^{(\sigma+\kappa)p}, 2^{(\sigma+\kappa+1)p})\}.$$

Then

$$\mu_d(\{r : \mathfrak{M}_\ell g > 2^\sigma\}) \lesssim \sum_{\kappa \ge 0} \sum_j \sum_{n \in \mathcal{B}_{\ell j \sigma}^\kappa} 2^{j(d-1)} 2^{n(d-1)} |D_n^{j+n}|.$$

Now

$$|D_n^{j+n}| \lesssim 2^j N(E^{j+n}, 2^j)$$

and

$$2^{jd} \lesssim 2^{-n(d-1)p} 2^{-(\sigma+\kappa)p} 2^{\ell(\frac{d-1}{p'}-\frac{1}{p})p} \|g\|_{L^p(I^{j+\ell})}^p \qquad \text{if } n \in \mathcal{B}_{\ell j \sigma}^\kappa.$$

Therefore by Minkowski's inequality

$$\Big(\sum_\sigma 2^{\sigma q} [\mu_d(\{r : \mathfrak{M}_\ell g > 2^\sigma\})]^{q/p}\Big)^{1/q}$$

$$\le \sum_{\kappa \ge 0} \Big(\sum_\sigma 2^{\sigma q} \Big[\sum_j \sum_{n \in \mathcal{B}_{\ell j \sigma}^\kappa} 2^{(j+n)(d-1)} |D_n^{j+n}|\Big]^{q/p}\Big)^{1/q}$$

$$\lesssim \sum_{\kappa \ge 0} \Big(\sum_\sigma 2^{\sigma q} \Big[\sum_j \sum_{n \in \mathcal{B}_{\ell j \sigma}^\kappa} 2^{\ell(\frac{d-1}{p'}-\frac{1}{p})p} 2^{-n(d-1)(p-1)} N(E^{j+n}, 2^j) 2^{-(\kappa+\sigma)p} \|g\|_{L^p(I^{j+\ell})}^p\Big]^{q/p}\Big)^{1/q}$$

$$\lesssim 2^{\ell(\frac{d-1}{p'}-\frac{1}{p})} \sum_{\kappa \ge 0} 2^{-\kappa} \Big(\sum_j \Big[\sum_\sigma \Big(\sum_{n \in \mathcal{B}_{\ell j \sigma}^\kappa} 2^{-n(d-1)(p-1)} N(E^{j+n}, 2^j)\Big)^{q/p}\Big]^{p/q} \|g\|_{L^p(I^{j+\ell})}^p\Big)^{1/p}. \tag{3.7}$$

Now observe that for every $\sigma$ there are at most three $n$ in $\mathcal{B}^\kappa_{\ell j \sigma}$ and for every $n$ there are at most three $\sigma$ such that $n \in \mathcal{B}^\kappa_{\ell j \sigma}$. Therefore

$$(3.8) \qquad \sum_\sigma \Big[\sum_{n \in \mathcal{B}^\kappa_{\ell j \sigma}} 2^{-n(d-1)(p-1)} N(E^{j+n}, 2^j)\Big]^{q/p} \lesssim \sup_j \sum_n 2^{-n(d-1)q/p'}[N(E^{j+n}, 2^j)]^{q/p} \lesssim A^q.$$

Consequently

$$\Big(\sum_\sigma 2^{\sigma q}[\mu_d(\{r : \mathfrak{M}_\ell g(r) > 2^\sigma\})]^{q/p}\Big)^{1/q} \lesssim A 2^{\ell(\frac{d-1}{p'} - \frac{1}{p})} \sum_{\kappa \geq 0} 2^{-\kappa} \Big(\sum_j \|g\|^p_{L^p(I^{j+\ell})}\Big)^{1/p}$$
$$\lesssim A 2^{\ell(\frac{d-1}{p'} - \frac{1}{p})} \|g\|_p.$$

We have shown that $\mathfrak{M}_\ell$ maps $L^p$ boundedly into $L^{pq}(\mu_d)$, with norm $O(2^{\ell(\frac{d-1}{p'} - \frac{1}{p})})$, for $1 \leq p < \infty$, $p \leq q < \infty$ (provided that $A < \infty$). The same argument applies to the case $q = \infty$, with only notational changes. Note that the operator norms $\|\mathfrak{M}_\ell\|$ are controlled by a geometric sequence converging to 0 if $p < d/(d-1)$. Since for $1 < p < \infty$ the Lorentz-spaces carry a Banach-space topology we may sum in $\ell$ and the proposition is proved in the case $p > 1$. However using a result by Stein and N.Weiss [8] on summing functions in weak-$L^1$ one can extend the argument to cover the case $p = 1$ as well. $\square$

## 4. Estimates in higher dimensions, cont.

We now give a proof of Theorems 1.3 and 1.4. For $g \in L^p(\mathbb{R}_+)$ and fixed $p$, $1 < p < \infty$, $\lambda > 0$ we define an operator $\mathfrak{N} = \mathfrak{N}_{p,\lambda}$ by

$$\mathfrak{N}g(r) = \sum_{k \in \mathbb{Z}} \sum_{n \geq 0} \chi_{D^k_n}(r) 2^{-k(\lambda-1)} \int_{2^{k-n}}^{2^{k+1}} s^{-1/p'} g(s)\, ds$$

and let $d\mu_\lambda = r^{\lambda-1} dr$.

**Proposition 4.1.** *Let $E \subset (0, \infty)$ such that $|\overline{E}| = 0$ and $1 < p < \infty$, $p \leq q \leq \infty$. Define $\Gamma_\lambda(\beta)$, $\mathfrak{v}_\lambda$ as in (2.8-9). Suppose that*

$$\sup_{|I| \geq 1} \Big(\int_0^\infty \beta^q \Big[\frac{\mathfrak{v}_\lambda(\Gamma_\lambda(\beta) \cap T(I))}{|I|}\Big]^{q/p} \frac{d\beta}{\beta}\Big)^{1/q} < \infty \qquad (\text{if } q < \infty)$$

$$\sup_{|I| \geq 1} \sup_\beta \beta \Big[\frac{\mathfrak{v}_\lambda(\Gamma_\lambda(\beta) \cap T(I))}{|I|}\Big]^{1/p} < \infty \qquad (\text{if } q = \infty).$$

*Then $\mathfrak{N}$ maps $L^p$ boundedly into $L^{pq}(\mu_\lambda)$.*

**Proof of Theorems 1.3 and 1.4.** In view of Corollary 2.6 the $L^p_{\text{rad}} \to L^{pq}$ boundedness of $\mathcal{M}_E$ is a direct consequence of Lemmas 3.2 and 3.3 and Proposition 4.1; the latter is applied for $p = p_d = d/(d-1)$, $\lambda = d$ and $g(s) = f_0(s) s^{(d-1)/p_d}$. $\square$

The special case $D = 1$ of the following result concerning averages turns out to be crucial in the proof of Proposition 4.1.

**Proposition 4.2.** *For $(x, t) \in \mathbb{R}^D \times \mathbb{R}_+$ define*

$$Sf(x, t) = u(x, t) \int_{|y-x| \leq 2t} f(y)\, dy$$

where $u$ is a nonnegative measurable function. Let
$$\Gamma(\beta) = \{(x,t) : t^D u(x,t) > \beta\}$$
and let $\mu$ be a positive measure in $\mathbb{R}^D \times \mathbb{R}_+$. Suppose $1 < p < \infty$, $p \le q \le \infty$ and that
$$\sup_Q \Big(\int_0^\infty \Big[\beta^p \frac{\mu(\Gamma(\beta) \cap T(Q))}{|Q|}\Big]^{q/p} \frac{d\beta}{\beta}\Big)^{1/q} < \infty \qquad \text{if } q < \infty,$$
$$\sup_Q \sup_\beta \beta \Big[\frac{\mu(\Gamma(\beta) \cap T(Q))}{|Q|}\Big]^{1/p} < \infty \qquad \text{if } q = \infty$$
holds; here we take the supremum over all cubes in $\mathbb{R}^D$ and $T(Q)$ is the cube in $\mathbb{R}^D \times \mathbb{R}_+$ with bottom $Q$. Then
$$\|Sf\|_{L^{pq}(\mathbb{R}^{D+1}_+, d\mu)} \le C\|f\|_{L^p(\mathbb{R}^D)}.$$

**Proof.** For $m \in \mathbb{Z}$ let
$$\Omega_m = \{x : Mf(x) > 2^m\}$$
where $Mf$ is the Hardy-Littlewood maximal function of $f$. Let $\{Q_\nu^m\}$ be a Whitney-decomposition of $\Omega_m$; here we assume that the Whitney cubes are binary cubes such that the coordinates of the corners are of the form $k_1 2^{k_2}$ with $k_1, k_2 \in \mathbb{Z}$. Define
$$\mathcal{R}_\nu^m = T(Q_\nu^m) \setminus \bigcup_{Q_{\nu'}^{m+1} \subset Q_\nu^m} T(Q_{\nu'}^{m+1}).$$
Then, if $f \ne 0$, it is easy to see that every $(x,t) \in \mathbb{R}^{D+1}_+$ belongs to some $T(Q_\nu^m)$ (for suitable $m$ depending on $x$ and $t$) and thus
$$\mathbb{R}^{D+1}_+ = \bigcup_{m,\nu} \mathcal{R}_\nu^m.$$
Let
$$E_\alpha = \{(x,t) : |Sf(x,t)| > \alpha\}$$
then it follows that
$$\mu(E_\alpha) = \sum_m \sum_\nu \mu(E_\alpha \cap \mathcal{R}_\nu^m).$$
If $(x,t) \in \mathcal{R}_\nu^m$ then we may pick $x_0$ such that $|x - x_0| \le c_1 t$ ($c_1$ is some geometrical constant) and such that $x_0 \notin \Omega^{m+1}$ which means $Mf(x_0) \le 2^{m+1}$. Therefore
$$t^{-D} \int_{|y-x| \le 2t} f(y)\,dy \le c_2 Mf(x_0) \le c_2 2^{m+1}.$$
Consequently if $|Sf(x,t)| > \alpha$ and $(x,t) \in \mathcal{R}_\nu^m$ then $u(x,t)t^D > \alpha(c_2 2^{m+1})^{-1}$ or $E_\alpha \cap \mathcal{R}_\nu^m \subset \Gamma(c_2^{-1} 2^{-m-1}\alpha)$. Thus
$$\mu(E_\alpha \cap \mathcal{R}_\nu^m) \le \mu(\Gamma(c_2^{-1} 2^{-m-1}\alpha) \cap T(Q_\nu^m))$$
and therefore
$$\Big(\int_0^\infty [\alpha^p \mu(E_\alpha)]^{q/p} \frac{d\alpha}{\alpha}\Big)^{1/q} \lesssim \Big(\int_0^\infty \Big[\alpha^p \sum_m \sum_\nu \mu(\Gamma(C 2^{-m}\alpha) \cap T(Q_\nu^m))\Big]^{q/p} \frac{d\alpha}{\alpha}\Big)^{1/q}$$
$$\lesssim \Big(\sum_m \sum_\nu \Big(\int_0^\infty \alpha^q [\mu(\Gamma(C 2^{-m}\alpha) \cap T(Q_\nu^m))]^{q/p} \frac{d\alpha}{\alpha}\Big)^{p/q}\Big)^{1/p}$$
$$= \Big(\sum_m \sum_\nu (C 2^m)^p |Q_\nu^m| \Big(\int_0^\infty \beta^q \Big[\frac{\mu(\Gamma(\beta) \cap T(Q_\nu^m))}{|Q_\nu^m|}\Big]^{q/p} \frac{d\beta}{\beta}\Big)^{p/q}\Big)^{1/p}$$
$$\lesssim \Big(\sum_m \sum_\nu 2^{mp} |Q_\nu^m|\Big)^{1/p}$$
$$= \Big(\sum_m 2^{mp} |\Omega^m|\Big)^{1/p} \lesssim \|Mf\|_p$$

and since $p > 1$ the asserted inequality follows from the $L^p$-boundedness of the Hardy-Littlewood maximal function. □

**Proof of Proposition 4.1.** We apply Proposition 4.2 with $D = 1$ and

$$f(x) = \sum_{\kappa \in \mathbb{Z}} \chi_{[\kappa,\kappa+1)}(x) \int_{2^\kappa}^{2^{\kappa+1}} s^{-1/p'} g(s) ds,$$

$$u(x,t) = \sum_k \chi_{[k,k+1)}(x) 2^{-k(\lambda-1)},$$

$$d\mu(x,t) = \sum_{k \in \mathbb{Z}} \sum_{n>0} \chi_{[k,k+1)}(x) \chi_{[n,n+1)}(t) |D_n^k| 2^{k(\lambda-1)} dx\, dt.$$

Let

$$U_\alpha = \{(k,n) : 2^{-k(\lambda-1)} \int_{2^{k-n}}^{2^k} s^{-1/p'} |g(s)| ds > \alpha\}.$$

Then

$$u(x,t) \int_{x-2t}^{x+2t} f(y) dy > \alpha \quad \text{if } (k,n) \in U_\alpha,\, k \leq x \leq k+1,\, n \leq t \leq n+1.$$

Therefore an application of Proposition 4.2 yields that under our hypothesis

$$\left(\int_0^\infty \alpha^q \Big[ \sum_{(k,n) \in U_\alpha} |D_n^k| 2^{k(\lambda-1)} \Big]^{q/p} \frac{d\alpha}{\alpha}\right)^{1/q} \lesssim \left(\sum_{\kappa \in \mathbb{Z}} \Big[\int_{2^\kappa}^{2^{\kappa+1}} s^{-1/p'} g(s) ds\Big]^p\right)^{1/p} \lesssim \|g\|_p$$

which implies the assertion. □

## 5. Estimates in two dimensions

Again we begin by stating a pointwise inequality for $\mathcal{M}_E$ acting on radial functions $f$ in $\mathbb{R}^2$ with $f(x) = f_0(r)$ where $r = |x|$.

**Lemma 5.1.** *Fix $1 \leq p < 2$ and set*

(5.1) $$g(s) = f_0(s) s^{1/p}.$$

*Then*

(5.2) $$\mathcal{M}_E f(x) \leq C[\mathfrak{M}g(r) + \widetilde{\mathfrak{M}}g(r) + \sum_{i=1}^4 R_i f_0(r)]$$

*where*

(5.3) $$\mathfrak{M}g(r) = \sup_{\substack{t \in E \\ r/2 < t < 3r/2}} r^{-1} \int_{|r-t|}^{r+t} s^{1/2-1/p}(s - |r-t|)^{-1/2} g(s) ds$$

(5.4) $$\widetilde{\mathfrak{M}}g(r) = \sup_{\substack{t \in E \\ r/2 < t < 3r/2}} r^{-1} \int_{|r-t|}^{r+t} s^{1/2-1/p}(r+t-s)^{-1/2} g(s) ds$$

and

$$(5.5) \quad R_1 f_0(r) = \sup_{\substack{t \in E \\ t \leq r/2}} t^{-1/2} \int_{r-t}^{r} |s - r + t|^{-1/2} f_0(s)\, ds$$

$$(5.6) \quad R_2 f_0(r) = \sup_{\substack{t \in E \\ t \leq r/2}} t^{-1/2} \int_{r}^{r+t} |r + t - s|^{-1/2} f_0(s)\, ds$$

$$(5.7) \quad R_3 f_0(r) = \sup_{\substack{t \in E \\ t \geq 3r/2}} r^{-1/2} \int_{t-r}^{t} |s - t + r|^{-1/2} f_0(s)\, ds$$

$$(5.8) \quad R_4 f_0(r) = \sup_{\substack{t \in E \\ t \geq 3r/2}} r^{-1/2} \int_{t}^{t+r} |r + t - s|^{-1/2} f_0(s)\, ds.$$

The proof consists of straightforward manipulations of (2.19-20) and is omitted. The case $d = 2$ of Theorems 1.1 and 1.2 follows from the results on $\mathfrak{M}, \widetilde{\mathfrak{M}}$ and $R_i$ stated in Propositions 5.2-5.4 below.

**Proposition 5.2.** *For $1 < p \leq \infty$ and $i = 1, 2$*

$$\left( \int_0^\infty |R_i f_0(r)|^p r\, dr \right)^{1/p} \leq C \sum_{m \geq 0} \sup_k [N(E^k, 2^{k-m})]^{1/p} 2^{-m/2} (m+1)^{1/p} \left( \int_0^\infty |f_0(r)|^p r\, dr \right)^{1/p};$$

*moreover there is the weak-type inequality*

$$\mu_2(\{r : R_i f_0(r) > \alpha\}) \leq C \sum_{m \geq 0} 2^{-m/2}(m+1) \sup_k N(E^k, 2^{k-m}) \alpha^{-1} \int_0^\infty |f_0(r)| r\, dr.$$

**Proof.** We only consider $R_1$; the operator $R_2$ is handled analogously.

We first observe that if $\operatorname{supp} h \subset [2^k, 2^{k+1}]$ then $\operatorname{supp} R_1 h \in [2^{k-2}, 2^{k+5}]$ and we may hence assume that $f_0$ is supported in $[2^k, 2^{k+1}]$. In this case

$$R_1 f_0(r) \lesssim \sup_{t \in E} K_t * f_0(r)$$

and the convolution kernel is defined by

$$K_t(x) = t^{-1}(1 - t^{-1}x)_+^{-1/2} \chi_{[0,\infty]}(x).$$

Let $\beta \in C_0^\infty$ be supported in $(1/2, 2)$ such that $\sum_{k=-\infty}^\infty \beta(2^k s) = 1$ and define for $m > 0$

$$K_t^m(x) = K_t(x) \beta(2^m(1 - t^{-1}x)).$$

Then

$$|R_1 f_0(r)| \lesssim M f_0(r) + \sum_{m=1}^\infty R_{1,m} f_0(r)$$

where $M$ is the Hardy-Littlewood maximal operator and $R_{1,m} f_0 = \sup_{t \in E} |K_t^m * f_0|$. Since $R_{1,m}$ is bounded on $L^\infty$ with operator norm $O(2^{-m/2})$ it suffices to prove the weak-type $(1,1)$ inequality with respect to Lebesgue measure. In view of the above mentioned properties of the support of $f_0$ and $R_1 f_0$ the analogous weighted version is an immediate consequence.

In order to prove the weak-type $(1,1)$ inequality we use Calderón-Zygmund theory for vector-valued operators. It then suffices to check that

$$(5.9) \qquad \int_{|x|>2|y|} \sup_k \sup_{t\in E^k} |K_t^m(x-y) - K_t^m(x)|\,dx \lesssim m 2^{-m/2} \sup_k N(E^k, 2^{k-m}).$$

First observe that $\operatorname{supp} K_t^m \subset \{x : |x| \leq t\}$ and therefore only terms with $2^k \geq |y|/8$ enter in the integral (5.9). Choose a minimal cover of $E^k$ with binary intervals $I_\nu^{km}$ of length $2^{k-m}$. It is elementary to check that

$$\sup_{t\in I_\nu^{km}} \left[ |K_t^m(x)| + 2^{k-m} \left|\frac{dK_t^m}{dx}(x)\right|\right] \leq C 2^{-k+m/2}$$

where $C$ does not depend on $k$; moreover for fixed $\nu$ the expression on the left hand side is supported in $J_\nu^{km} = \{x \in \mathbb{R} : \operatorname{dist}(x, I_\nu^{km}) \leq 2^{k+1-m}\}$ which is an interval of length $\leq 2^{k-m+3}$. Therefore the left hand side of (5.9) is bounded by

$$2 \sum_{|y|/8 \leq 2^k \leq 2^m|y|} \sum_\nu \int_{J_\nu^{km}} \sup_{t\in I_\nu^{km}} |K_t^m(x)|\,dx + \sum_{2^k > 2^m|y|} \sum_\nu \int_0^1 \int_{J_\nu^{km}} \sup_{t\in I_\nu^{km}} \left|y\frac{dK_t^m}{dx}(x-sy)\right|\,dx\,ds$$

$$\lesssim \sum_{|y|/8 \leq 2^k \leq 2^m|y|} \sum_\nu 2^{-k+m/2} |J_\nu^{km}| + \sum_{2^k > 2^m|y|} |y|2^{m-k} \sum_\nu |J_\nu^{km}| 2^{-k+m/2}$$

$$\lesssim m 2^{-m/2} \sup_k N(E^k, 2^{k-m}). \quad \square$$

The following result concerning $R_3$ and $R_4$ is a singular variant of Lemma 3.3.

**Proposition 5.3.** *For $1 \leq p < \infty$ and $i = 3, 4$*

$$\left(\int_0^\infty |R_i f_0(r)|^p r\,dr\right)^{1/p} \leq C \sum_{m\geq 0} \sup_k [N(E^k, 2^{k-m})]^{1/p} 2^{-m/2} \left(\int_0^\infty |f_0(r)|^p r\,dr\right)^{1/p}.$$

**Proof.** We only consider $R_3$; the corresponding proof for $R_4$ is similar. Define

$$R_{3,m} f_0(r) = \sup_{\substack{t\in E \\ t\geq 3r/2}} 2^{m/2} r^{-1} \int_{t-r+2^{-m-1}r}^{t-r+2^{-m}r} |f_0(s)|\,ds.$$

Then $|R_3 f_0(r)| \lesssim \sum_{m=0}^\infty |R_{3,m} f_0(r)|$. We use Lemma 3.3 to see that the operators $R_{3,m}$ are bounded for $m \leq 4$ and assume henceforth $m > 4$. For fixed $m > 4$ we introduce a further decomposition in terms of the $t/r$; we then have

$$R_{3,m} f_0(r) \lesssim \sum_{L=0}^\infty R_{3,m,L} f_0(r)$$

where

$$R_{3,m,L} f_0(r) = \sum_{k\in\mathbb{Z}} \chi_{I_k}(r) 2^{-2k} 2^{m/2} 2^{-L} \sup_{t\in E^{k+L}} \int_{t-r+2^{-m-1}r}^{t-r+2^{-m}r} |f_0(s)| s\,ds.$$

Note that the operator norm of $R_{3,m,L}$ on $L^\infty$ is bounded by $C 2^{-m/2}$ uniformly in $L$. We shall prove that $R_{3,m,L}$ is bounded on $L^1(\mu_2)$ with operator norm bounded by $C 2^{-L} 2^{-m/2} \sup_k N(E^k, 2^{k-m-L})$. Taking this for granted we obtain by interpolation that the $L^p(\mu_2)$ operator norm of $R_{3,m,L}$ is bounded by

$C2^{-L/p}2^{-m/2}\sup_k[N(E^k,2^{k-m-L})]^{1/p}$ and this implies that

$$\left(\int_0^\infty |R_3 f_0(r)|^p r\,dr\right)^{1/p} \lesssim \sum_{L>0} 2^{-L/p}\sum_{m\geq 0}\sup_k[N(E^k,2^{k-m-L})]^{1/p}2^{-m/2}\left(\int_0^\infty |f_0(r)|^p r\,dr\right)^{1/p}$$

$$\lesssim \sum_{n\geq 0}\sup_k[N(E^k,2^{k-n})]^{1/p}2^{-n/2}\left(\int_0^\infty |f_0(r)|^p r\,dr\right)^{1/p}.$$

In order to prove the required $L^1(\mu_2)$ inequality for $R_{3,m,L}$ we cover the set $E^{k+L}$ with intervals $I_\nu^{k+L,m+L}$ of length $2^{k-m}=2^{k+L-(m+L)}$. Denote by $\widetilde{I}_\nu^{k+L,m+L}$ the double interval. Then

$$\int_{2^k}^{2^{k+1}}|R_{3,m,L}f_0(r)|r\,dr \lesssim 2^{-2k}2^{m/2}2^{-L}\sum_\nu \int_{2^k}^{2^{k+1}}\int_{-r+\widetilde{I}_\nu^{k+L,m+L}}|f_0(s)|s\,ds\,r\,dr$$

$$\lesssim \sum_\nu |I_\nu^{k+L,m+L}|2^{-k-L+m/2}\int_{2^{k+L}-1}^{2^{k+L+3}}|f_0(s)|s\,ds$$

$$\lesssim 2^{-L-m/2}N(E^{k+L},2^{k+L-(m+L)})\int_{2^{k+L}-1}^{2^{k+L+3}}|f_0(s)|s\,ds.$$

Summing in $k$ yields the asserted $L^1$ inequality. $\square$

**Proposition 5.4.** *Suppose that*

$$\sup_j\left(\sum_{n\geq 0}[N(E^{j+n},2^j)^{1/p}2^{-n/p'}]^q\right)^{1/q} \leq A \qquad \text{if } p\leq q<\infty,$$

$$\sup_{j,n} N(E^{j+n},2^j)^{1/p}2^{-n/p'} \leq A \qquad \text{if } q=\infty.$$

*Then for $1\leq p<2$, $p\leq q\leq\infty$*

$$\|\mathfrak{M}g\|_{L^{pq}(\mu_d)} + \|\widetilde{\mathfrak{M}}g\|_{L^{pq}(\mu_d)} \lesssim A\|g\|_{L^p}.$$

**Proof.** We show this only for the operator $\mathfrak{M}$; the proof for $\widetilde{\mathfrak{M}}$ is similar but simpler. We dominate

$$\mathfrak{M}g(r) \lesssim \mathfrak{M}_0 g(r) + \sum_{\ell=1}^\infty \mathfrak{M}_\ell g(r)$$

where

$$\mathfrak{M}_0 g(r) = \sup_{\substack{t\in E \\ r/2<t<3r/2}} r^{-1}\int_{|r-t|}^{2|r-t|} s^{1/2-1/p}(s-|r-t|)^{-1/2}g(s)ds$$

and for $\ell>0$ the operators $\mathfrak{M}_\ell$ are defined in (3.6). We see that for $\ell\geq 1$ the estimate

$$\|\mathfrak{M}_\ell g\|_{L^{pq}(\mu_2)} \lesssim 2^{\ell(1/p'-1/p)}A\|g\|_p$$

was already obtained in the proof of Proposition 3.4. However the term $\mathfrak{M}_0 g$ is more singular in two dimensions and in what follows we prove the required estimate for the operator $\mathfrak{M}_0$.

The set $D_n^k$ is a union of intervals $I_{n\nu}^k$. Let for $0\leq l\leq n$

$$E_{nl\nu}^k = \{t\in E^k : 2^{k-n+l}\leq \text{dist}(t,I_{n\nu}^k)\leq 2^{k-n+l+1}\}.$$

Then a straightforward estimate yields

$$\mathfrak{M}_0 g(r)$$
$$\leq C_1 \sum_{k\in\mathbb{Z}} \sum_{n\geq 0} \sum_\nu \chi_{I^k_{n\nu}}(r) \sum_{l=0}^n \sum_{m=0}^\infty 2^{-2k/p} 2^{(2n-2l+m)(1/p-1/2)} \sup_{t\in E^k_{nl\nu}} \left( \int_{|r-t|+2^{k-n+l-m}}^{|r-t|+2^{k-n+l-m+1}} |g(s)|^p ds \right)^{1/p}$$
$$=: C_1 \sum_{l,m=0}^\infty \mathfrak{M}_{lm} g(r)$$

where we write

$$\mathfrak{M}_{lm} g(r) = \sum_{n\geq l} \sum_k \sum_\nu \chi_{I^k_{n\nu}}(r) 2^{-2k/p} 2^{(2n-2l+m)(1/p-1/2)} \sup_{t\in E^k_{nl\nu}} \left( \int_{|r-t|+2^{k-n+l-m}}^{|r-t|+2^{k-n+l-m+1}} |g(s)|^p ds \right)^{1/p}$$
$$= \sum_j \sum_{n\geq l} \sum_\nu \chi_{I^{j+n}_{n\nu}}(r) 2^{-2(j+n)/p} 2^{(2n-2l+m)(1/p-1/2)} \sup_{t\in E^{j+n}_{nl\nu}} \left( \int_{|r-t|+2^{j+l-m}}^{|r-t|+2^{j+l-m+1}} |g(s)|^p ds \right)^{1/p}.$$

Let $I_{j+l} = [2^{j+l}, 2^{j+l+1}]$ and denote by $\widetilde{I}_{j+l}$ the expanded interval as defined in the introduction. We proceed analogously as in the proof of Proposition 3.4 and define

$$\mathcal{B}^{\kappa lm}_{j\sigma} = \{n \geq 0 : C_1 2^{-2(j+n)/p} 2^{(2n-2l+m)(1/p-1/2)} \|g\|_{L^p(\widetilde{I}_{j+l})} \in [2^{\sigma+\kappa}, 2^{\sigma+\kappa+1}]\}$$

There is the crude estimate

(5.10) $$\mathfrak{M}_{lm} g(r) \leq 2^{-2(j+n)/p} 2^{(2n-2l+m)(1/p-1/2)} \|g\|_{L^p(\widetilde{I}_{j+l})}, \qquad r \in D^{j+n}_n,$$

which implies that

$$\mu_2(\{r : \mathfrak{M}_{lm} g(r) > 2^\sigma\}) \lesssim \sum_{\kappa \geq 0} \sum_j \sum_{n\in\mathcal{B}^{\kappa lm}_{j\sigma}} 2^{j+n} |D^{j+n}_n|.$$

Here we have of course used that $|\overline{E}| = 0$. Now

$$\left( \sum_\sigma 2^{\sigma q} [\mu_2(\{r : \mathfrak{M}_{lm} g(r) > 2^\sigma\})]^{q/p} \right)^{1/q}$$
$$\leq \sum_{\kappa \geq \max\{m-l,0\}} \left( \sum_\sigma 2^{\sigma q} \Big[\sum_j \sum_{n\in\mathcal{B}^{\kappa lm}_{j\sigma}} 2^{j+n} |D^{j+n}_n|\Big]^{q/p} \right)^{1/q}$$
$$+ \sum_{\kappa < \max\{m-l,0\}} \left( \sum_\sigma 2^{\sigma q} \Big[\sum_j \sum_{n\in\mathcal{B}^{\kappa lm}_{j\sigma}} 2^{j+n} |\{r \in D^{j+n}_n : \mathfrak{M}_{lm} g(r) > 2^\sigma\}|\Big]^{q/p} \right)^{1/q}$$
$$= \mathcal{E}^1_{lm} + \mathcal{E}^2_{lm}$$

where the expression $\mathcal{E}^2_{lm}$ is of course only present when $m > l$. The estimate of $\mathcal{E}^1_{lm}$ follows the lines of the proof of Proposition 3.4 while for $\mathcal{E}^2_{lm}$ we shall use the finer estimate

(5.11) $$\left( 2^k \int_{D^k_n} |\mathfrak{M}_{lm} g(r)|^p dr \right)^{1/p} \leq C 2^{(l-n)/p'} 2^{-m/2} [N(E^k, 2^{k-n+l-m-8})]^{1/p} \|g\|_{L^p(\widetilde{I}_{k-n+l})}.$$

Now if $n \in \mathcal{B}^{\kappa lm}_{j\sigma}$ then

$$2^{2j+n} \leq C_1^p (2^{-\kappa-\sigma})^p 2^{-n(p-1)} 2^{l(p-2)} 2^{m(1-p/2)} \|g\|^p_{L^p(\widetilde{I}_{j+l})}$$

and therefore we obtain by arguing as in (3.7), (3.8)

$$\mathcal{E}_{lm}^1 \lesssim \sum_{\kappa \geq \max(m-l,0)} 2^{-\kappa} 2^{(m-2l)(1/p-1/2)} \Big(\sum_\sigma \Big[\sum_j \sum_{n \in \mathcal{B}_{j\sigma}^{\kappa lm}} 2^{-n(p-1)} N(E^{j+n}, 2^j) \|g\|_{L^p(\widetilde{I}_{j+l})}^p\Big]^{q/p}\Big)^{1/q}$$
(5.12)
$$\lesssim \min\{2^{(m-2l)(1/p-1/2)}, 2^{m(1/p-3/2)} 2^{2l(1-1/p)}\} \sup_j \Big(\sum_n [2^{-n(p-1)} N(E^{j+n}, 2^j)]^{q/p}\Big)^{1/q} \|g\|_p.$$

Assuming (5.11) we obtain for $m > l$ by Tshebyshev's inequality

$$\mathcal{E}_{lm}^2 \leq \sum_{\kappa < m-l} \Big(\sum_\sigma \Big[\sum_j \sum_{n \in \mathcal{B}_{j\sigma}^{\kappa lm}} 2^{(l-n)(p-1)} 2^{-mp/2} N(E^{j+n}, 2^{j+l-m-8}) \|g\|_{L^p(\widetilde{I}_{j+l})}^p\Big]^{q/p}\Big)^{1/q}$$

$$\lesssim (m-l) \sup_j \Big(\sum_n [2^{(l-n)(p-1)} 2^{-mp/2} N(E^{j+n}, 2^{j+l-m-8})]^{q/p}\Big)^{1/q} \|g\|_p$$

(5.13)
$$\lesssim (m-l) 2^{-m(1/p-1/2)} \sup_j \Big(\sum_n [2^{-n(p-1)} N(E^{j+n}, 2^j)]^{q/p}\Big)^{1/q} \|g\|_p.$$

¿From (5.12) and (5.13) it follows that

$$\sum_{m \geq 0} \sum_{l \geq 0} \mathcal{E}_{lm}^1 + \sum_{m \geq 0} \sum_{l < m} \mathcal{E}_{lm}^2 \lesssim \sum_m (1+m)^2 2^{-m(1/p-1/2)} \sup_j \Big(\sum_n [2^{-n(p-1)} N(E^{j+n}, 2^j)]^{q/p}\Big)^{1/q} \|g\|_p$$

which proves the asserted inequality for the case $1 < p < 2$, $p \leq q < \infty$. The argument for $q = \infty$ is analogous and the case $p = 1$ can be handled in a similar way using the result of [8], as in the proof of Proposition 3.4.

*Proof of (5.11).* For small $m$ the appropriate estimates have already been obtained in the proof of Proposition 3.4. Therefore we may assume that $m \geq 4$ in what follows.

For an integer $\mu$ with $2^{n-l+4} \leq \mu \leq 2^{n-l+5}$ let $J_\mu^{knl}$ be the union of all intervals $I_{n\nu}^k$ (in $D_n^k$) which have nonempty intersection with $[\mu 2^{k-n+l-4}, (\mu+1) 2^{k-n+l-4}]$ and let

$$F_\mu^{knl} = \{t \in E^k : 2^{k-n+l-1} \leq \mathrm{dist}(t, J_\mu^{knl}) \leq 2^{k-n+l+1}\}.$$

Thus $E_{nl\nu}^k \subset F_\mu^{knl}$ whenever $I_{n\nu}^k$ has nonempty intersection with $[\mu 2^{k-n+l-4}, (\mu+1) 2^{k-n+l-4}]$.

We choose a minimal cover of $F_\mu^{knl}$ with binary intervals $Q_{\mu\rho}^{knl} = [a_\rho, b_\rho]$ of length $2^{k-n+l-m-8}$. If $r \in J_\mu^{knl}$ then for given $\rho$ we have either $b_\rho < r$ or $r < a_\rho$. Define for $r \in J_\mu^{knl}$

$$\widetilde{Q}_{\mu\rho}^{knl}(r) = \begin{cases} [r - b_\rho + 2^{k-n+l-m-3}, r - a_\rho + 2^{k-n+l-m-2}] & \text{if } b_\rho < r \\ [a_\rho - r + 2^{k-n+l-m-3}, b_\rho - r + 2^{k-n+l-m-2}] & \text{if } r < a_\rho \end{cases}.$$

For fixed $s$ let

$$R_{\mu\rho}^{knl}(s) = \{r \in J_\mu^{knl} : s \in \widetilde{Q}_{\mu\rho}^{knl}(r)\}.$$

Then the measure of $R_{\mu\rho}^{knl}(s)$ is $O(2^{k-n+l-m})$. We estimate

$$\Big(2^k \int_{D_n^k} |\mathfrak{M}_{lm} g(r)|^p dr\Big)^{1/p}$$

(5.14)
$$\leq C 2^{-2k/p} 2^{(2n-2l+m)(1/p-1/2)} \Big(2^k \sum_\mu \int_{J_\mu^{knl}} \sup_{t \in F_\mu^{knl}} \int_{|r-t|+2^{k-n+l-m}}^{|r-t|+2^{k-n+l-m+1}} |g(s)|^p \, ds \, dr\Big)^{1/p}$$

and

$$\sum_\mu \int_{J_\mu^{knl}} \sup_{t \in F_\mu^{knl}} \int_{|r-t|+2^{k-n+l-m}}^{|r-t|+2^{k-n+l-m+1}} |g(s)|^p \, ds \, dr$$

$$\lesssim \sum_\mu \sum_\rho \int_{J_\mu^{knl}} \int_{\widetilde{Q}_{\mu\rho}^{knl}(r)} |g(s)|^p \, ds \, dr$$

$$\lesssim \sum_\mu \sum_\rho \int_{2^{k-n+l}}^{2^{k-n+l+2}} \int_{R_{\mu\rho}^{knl}(s)} dr \, |g(s)|^p \, ds$$

$$\lesssim 2^{k-n+l-m} \sum_\mu N(F_\mu^{knl}, 2^{k-n+l-m-8}) \int_{2^{k-n+l}}^{2^{k-n+l+2}} |g(s)|^p \, ds$$

(5.15) $$\lesssim 2^{k-n+l-m} N(E^k, 2^{k-n+l-m-8}) \int_{2^{k-n+l}}^{2^{k-n+l+2}} |g(s)|^p \, ds.$$

Combining (5.14) and (5.15) we immediately get the desired estimate (5.11). □

## References

1. J. Bourgain, *Averages in the plane over convex curves and maximal operators*, Jour. Anal. **47** (1986), 69–85.
2. M. Leckband, *A note on the spherical maximal operator for radial functions*, Proc. Amer. Math. Soc. **100** (1987), 635–640.
3. A. Seeger, S. Wainger and J. Wright, *Pointwise convergence of spherical means*, Math. Proc. Camb. Phil. Soc. **118** (1995), 115–124.
4. E. M. Stein, *On limits of sequences of maximal operators*, Annals of Math. **74** (1961), 140–170.
5. \_\_\_\_\_\_, *Maximal functions: spherical means*, Proc. Nat. Acad. Sci. **73** (1976), 2174–2175.
6. \_\_\_\_\_\_, *Harmonic analysis: Real variable methods, orthogonality and oscillatory integrals*, Princeton Univ. Press, Princeton, 1993.
7. E. M. Stein and G. Weiss, *Introduction to Fourier analysis on Euclidean spaces*, Princeton Univ. Press, Princeton, N.J., 1971.
8. E. M. Stein and N. J. Weiss, *On the convergence of Poisson integrals*, Trans. Amer. Math. Soc. **140** (1969), 34–54.

DEPARTMENT OF MATHEMATICS, UNIVERSITY OF WISCONSIN, MADISON, WI 53706, USA (ADDRESS OF THE FIRST AND THE SECOND AUTHOR)

SCHOOL OF MATHEMATICS, UNIVERSITY OF NEW SOUTH WALES, KENSINGTON NSW 2033, AUSTRALIA (ADDRESS OF THE THIRD AUTHOR)